\documentclass[12pt,reqno]{amsart}

\usepackage{fullpage,amsfonts,amsmath,amssymb,amscd,youngtab}

\input xypic

\theoremstyle{plain}
\newtheorem*{lem}{Lemma}

\newtheorem*{assumption}{Assumption}
\newtheorem*{cor}{Corollary}
\newtheorem*{thrm}{Theorem}

\theoremstyle{Definition}

\newtheorem*{rem}{Remark}

\newcommand{\triv}{\mathsf{triv}}

\newcommand{\res}{\mathrm{res}}

\newcommand{\End}{\mathrm{End}}
\newcommand{\cO}{\mathcal{O}}
\newcommand{\irrep}{\textsf{Irrep}}

\newcommand{\rank}{\mathrm{rank}}

\newcommand{\e}{\varepsilon}

\newcommand{\C}{\mathbb{C}}

\newcommand{\KZ}{\texttt{KZ}}

\newcommand{\h}{\mathfrak{h}}

\title{Rational Cherednik algebras and diagonal coinvariants
 of $G(m,p,n)$}

\author{Richard Vale}

\address{Department of Mathematics, University
of Glasgow, Glasgow, G12 8QW, U.K.}

\email{rv@maths.gla.ac.uk}

\date{\today}

\begin{document}
\begin{abstract}
We construct a quotient ring of the ring of diagonal coinvariants
of the complex reflection group $W=G(m,p,n)$ and determine its
graded character. This generalises a result of Gordon for Coxeter
groups. The proof uses a study of category $\cO$ for the rational
Cherednik algebra of $W$.
\end{abstract}


 \maketitle

\section{Introduction}

\subsection{}
\label{gordon} Let $\mathfrak{h}$ be a finite--dimensional complex
vector space. An element $s \in \End (\h)$ is called a complex
reflection if $\rank_\h (1-s) =1$ and $s$ has finite order. A
finite group generated by complex reflections is called a complex
reflection group. If $W$ is a complex reflection group then the
ring of invariants $\C [ \h ]^W$ is a polynomial ring by the
Shepherd-Todd theorem \cite[Theorem 7.2.1]{benson} and if $\C
[\h]^{W}_{+}$ denotes the elements with zero constant term then it
is well-known that the ring of coinvariants
$$\frac{\C [\h]}{\langle \C [\h]^{W}_{+} \rangle}$$
is a finite--dimensional algebra isomorphic to $\C W$ as a
$W$--module. There is interest in analogues of this construction
with the representation $\h \oplus \h^*$ in place of $\h$, see for
example \cite{Hai}. The ring
$$D_W := \frac{\C [\h \oplus \h^*]}{\langle \C [\h \oplus \h^*]^{W}_{+} \rangle}$$
is called the ring of diagonal coinvariants of $W$.
The ring $D_W$ has a natural grading with $\deg(\h^*) =1$ and
$\deg(\h) = -1$. The following result was conjectured by Haiman
and proved in Gordon \cite{Go}:

\begin{thrm}\label{Gordon}\cite{Go}
Let $W$ be a finite Coxeter group of rank $n$ with Coxeter number
$h$ and sign representation $\e$. Then there exists a $W$--stable
quotient ring $R_W$ of $D_W$ with the properties:
\begin{enumerate}
\item $\dim(R_W) = (h+1)^n$. \item $R_W$ is graded with Hilbert
series $t^{-hn/2}(1+t + \cdots + t^h)^n$. \item The image of $\C
[\h]$ in $R_W$ is $\C [\h] / \langle \C [\h]^{W}_{+} \rangle$.
\item The character $\chi$ of the $W$--module $R_W \otimes \e $
satisfies $\chi(w) = (h+1)^{\dim \ker (1-w)} \quad \forall w \in
W$.
\end{enumerate}
\end{thrm}

\subsection{}
\label{mainthm} In \cite{vale}, this result was generalised to the
complex reflection groups $G(m,1,n)$. The aim of this paper is to
obtain a further generalisation to the groups $G(m,p,n)$, with
some mild restrictions on $m,p,n$. The following result will be
proved:

\begin{thrm}\label{mainthm}
Let $W = G(m,p,n)$ where $m \neq p$
 and let $\h$ be the reflection representation of $W$. Let $d=m/p$.
Then there exists a $W$--stable quotient ring $S_W$ of $D_W$ with
the properties:
\begin{enumerate}
\item $\dim(S_W) = (m(n-1)+d+1)^n$. \item $S_W$ is graded with
Hilbert series $t^{-n-m \binom{n}{2} }(1 + t + \cdots +
t^{m(n-1)+d})^n$. \item The image of $\C [\h]$ in $S_W$ is $\C
[\h] / \langle \C [\h]^{W}_{+} \rangle$. \item The character
$\chi$ of $S_W \otimes \wedge^n \h^*$ as a $W$--module satisfies
$\chi(w) = (m(n-1)+d+1)^{\dim \ker (1-w)} \quad \forall w \in W$.
\end{enumerate}
\end{thrm}

\subsection{}
Theorem \ref{gordon} is proved by obtaining $R_W$ as the
associated graded module of a finite--dimensional module over the
rational Cherednik algebra of $W$. The properties of this module
are derived by studying the category $\mathcal{O}$ for the
rational Cherednik algebra. This is also the method that will be
used to prove Theorem \ref{mainthm}.

\subsection{}
The structure of the paper is as follows. In Section
\ref{cherednik}, the rational Cherednik algebra $H_\kappa$ is
introduced for $W=G(m,p,n)$. This is a certain deformation of the
skew group algebra $\C[\h \oplus \h^*]*W$ which depends on
parameters $\kappa \in \C^{m/p}$. Next, in Section \ref{cato}, we
recall some important properties of category $\cO$ for $H_\kappa$,
including the Knizhnik-Zamolodchikov functor $\KZ$, which we use
to relate category $\cO$ to the category of modules over a Hecke
algebra $\mathcal{H}$. A parametrisation of the simple
$\mathcal{H}$--modules given by Genet and Ja\c{c}on \cite{GJ}
enables us to prove that, for suitable choices of the parameters,
there is only one finite-dimensional simple object
$L(\mathsf{triv})$ in category $\cO$ (Theorem \ref{infdim}). We
then define, in Section \ref{shifting}, a one-dimensional
$H_{\kappa^{'}}$--module $\Lambda^\psi$ for ``shifted" values of
the parameters $\kappa^{'}$, and, using a shift isomorphism due to
Berest and Chalykh \cite{BC}, we construct the shifted module
$$L:= H_{\kappa} e_\e \otimes_{e_\e H_{\kappa} e_\e} \Lambda^\psi$$
where $e_\e$ is a certain idempotent in $\C W$. The associated
graded module $\mathrm{gr} L$ is, up to tensoring by a
one-dimensional $W$--module, naturally a quotient of the ring of
diagonal coinvariants. But $L$ is also a finite-dimensional object
of category $\cO$, and we are able to use our results on category
$\cO$ to show that it is isomorphic to $L(\mathsf{triv})$. By
results of Chmutova and Etingof \cite{CE}, $L(\mathsf{triv})$ is
well-understood, and this enables us to compute the Hilbert series
and character of $L$ and hence of $\mathrm{gr} L$, proving Theorem
\ref{mainthm}. These calculations are given in Section
\ref{mainproof}.

\subsection{}
The main difference between our proof and the proof of Theorem
\ref{Gordon} is that for some of the Coxeter groups considered in
\cite{Go}, the rational Cherednik algebra depends on only one
parameter, and there is only one choice of parameter for which the
proof will work. In the $G(m,p,n)$ case with $m /neq p$, there is
greater freedom to choose the parameters, and hence it is possible
to have a lot of control over the simple modules in category
$\cO$, by considering what happens when the parameters are chosen
generically. We note here that in the $m=p$ case the rational
Cherednik algebra usually depends on only one parameter, and hence
the proof of Theorem \ref{mainthm} will not work for the groups
$G(m,m,n)$. However, it is likely that an analogue of Theorem
\ref{mainthm} can be proved for these groups, following the
arguments of \cite{Go}. We hope to return to this in future work.

\subsection{}
It may appear that, in the case $m=1$, the Hilbert series of $S_W$
should be $t^{-n-{n \choose 2}} (1 + t + \cdots + t^n )^n$, which
does not generalise Theorem \ref{Gordon} in type $A$. However, in
order to make Theorem \ref{mainthm} and Theorem \ref{Gordon} agree
in this case, we should write the Hilbert series of $S_W$ as
$t^{-\dim \h - m {\dim \h \choose 2}} (1 + t+ \cdots +
t^{m(n-1)+d})^{\dim \h}$. Note that the type $A$ case of Theorem
\ref{Gordon} does \emph{not} follow from Theorem \ref{mainthm}
because we will assume throughout that $m>1$.


\subsection{Acknowledgements}
 The research described here will form
part of the author's PhD thesis at the University of Glasgow. The
author thanks K. A. Brown and I. Gordon for suggesting this
problem and for their advice and encouragement. The author also
wishes to thank Y. Berest and O. Chalykh for allowing us to look
at a preliminary version of their paper \cite{BC}, and O. Chalykh
for explaining in detail to us the proof of Theorem \ref{shift}.

\section{The group $G(m,p,n)$}
\subsection{}
Let $\h$ be an $n$--dimensional complex vector space equipped with
a sesquilinear form $\langle -,- \rangle$. Let $m \ge 1$ and let
$p$ be a natural number such that $p|m$. We fix the notation
$d=m/p$ and $\e = e^{\frac{2 \pi i}{m}}$ throughout. The complex
reflection group $G(m,p,n)$ is defined to be the subgroup of
$GL(\h)$ consisting of those matrices with exactly one nonzero
entry in each row and column, such that the nonzero entries are
powers of $\e$ and such that the $d^{\mathrm{th}}$ power of the
product of the nonzero entries is $1$.

\subsection{Complex reflections}\label{def} We wish to identify the complex
reflections in $G(m,p,n)$. It turns out that the set of complex
reflections depends on $(m,p,n)$, so we will make the following
assumption
\begin{assumption}
From now on, assume $m>p$. Equivalently, $d>1$.
\end{assumption}
Under this assumption, the complex reflections in $G(m,p,n)$ are
as follows. If $\{y_1 , \ldots , y_n \}$ is an orthonormal basis
of $\h$ then  the set $\mathcal{S}$ of complex reflections in $W$
consists of the elements $s_i^{qp}$ for $1 \le i \le n$ and $1 \le
q \le d-1$, and $\sigma_{ij}^{(\ell)}$ for $1 \le i < j \le n$ and
$0 \le \ell \le m-1$, defined by:
\begin{eqnarray*}
 s_{i}^{qp} ( y_i) = \e^{qp} y_i & s_{i}^{qp} (y_j) = y_j & j \neq i \\
\sigma_{ij}^{(\ell)} (y_i) = \e^{-\ell} y_j & \sigma_{ij}^{(\ell)}
(y_j) =
\e^\ell y_i &  \sigma_{ij}^{(\ell)} (y_k) = y_k  , \: k \neq i,j .\\
\end{eqnarray*}

\subsection{}\label{def2}
We now list the $W$--conjugacy classes in $\mathcal{S}$. For each
$q$, $\{ s_i^{qp} | 1 \le i \le n\}$ form a single conjugacy class
in $\mathcal{S}$. If $n \ge 3$ or $n=2$ and $p$ is odd, then $\{
\sigma_{ij}^{(\ell)} | i < j, \: 0 \le \ell \le m-1 \}$ also form
a single conjugacy class in $\mathcal{S}$. For convenience, we now
make the following assumption.
\begin{assumption}
Either $n \ge 3$, or $n=2$ and $p$ is odd.
\end{assumption}
Under assumptions \ref{def} and \ref{def2}, we see that there are
exactly $d$ $W$--conjugacy classes of complex reflections.
Furthermore, it follows from \cite[Section 3]{BMR} that the
defining representation $\h$ is irreducible when these assumptions
hold. We are now in a position to construct the rational Cherednik
algebra of $W$.
\begin{rem}
Theorem \ref{mainthm} still holds in the case where $n=2$ and $p$
is even. To avoid clutter, the modifications necessary to prove
this case are explained in Section \ref{fudge}.
\end{rem}

\section{The rational Cherednik algebra}\label{cherednik}
In \cite{DO} and \cite{GGOR} we have the following definition. Let
$\mathcal{A}$ be the set of reflection hyperplanes of $W$ and for
$H \in \mathcal{A}$ let $W_H = \mathrm{stab}_W(H)$, a cyclic group
of order $e_H$. For $1 \le i \le e_H-1$, let $\e_{H,i} =
\frac{1}{e_H} \sum_{w \in W_H} \det (w)^i w$. For $H \in
\mathcal{A}$, let $\{k_{H, i} \}_{i=0}^{e_H}$ be a family of
scalars such that $k_{H,i} = k_{H',i}$ whenever $H, H'$ are in the
same $W$--orbit, and $k_{H,0} = k_{H, e_H} = 0$ for all $H$. For
each $H \in \mathcal{A}$, pick a linear form $\alpha_H \in \h^*$
with kernel $H$, and choose $\alpha_H^{\vee} \in \h$ such that $\C
\alpha_H^{\vee}$ is a $W_H$--stable complement to $H$ and
$\alpha_H (\alpha_H^{\vee}) = 2$.

The rational Cherednik algebra is defined to be the quotient of
$T(\h \oplus \h^*) \ast W$, the skew product of $W$ with the
tensor algebra on $\h \oplus \h^*$, by the relations $[y_1,
y_2]=0$ for all $y_1, y_2 \in \h$, $[x_1, x_2]=0$ for all $x_1,
x_2 \in \h^*$, and \begin{equation}\label{commute} [y,x]= \langle
y, x \rangle + \sum_{H \in \mathcal{A}} \frac{\langle \alpha_H, y
\rangle \langle \alpha_H^{\vee}, x \rangle}{2} e_H
\sum_{j=0}^{e_H-1} (k_{H, j+1} - k_{H,j} ) \e_{H,j}\end{equation}
for all $y \in \h$ and all $x \in \h^*$, where $\langle -,-
\rangle$ here denotes the evaluation pairing between $\h$ and
$\h^*$.

\begin{rem}
In this paper, the signs in the commutation relation have been
chosen so that our parameters $k_{H,i}$ are the same as those of
the paper \cite{GGOR}. However, we will also use results from the
paper \cite{BC}, in which the parameters denoted $k_{H,i}$ are the
negatives of those given here. See Section \ref{shift}.
\end{rem}

\subsection{}\label{commutationrel}
In the case $W=G(m,p,n)$, we may write out the commutation
relation more explicitly. Let $H_i$ be the reflection hyperplane
of $s_i^{qp}$ and $H_{ij \ell}$ be the reflection hyperplane of
$\sigma_{ij}^{(\ell)}$. Let $\{ y_1, \ldots , y_n\}$ be the
standard basis of $\h$ and $\{x_1, \ldots , x_n\}$ the dual basis
of $\h^*$. Then we may choose $\alpha_{H_i} = x_i$,
$\alpha_{H_i}^{\vee} = 2y_i$, $\alpha_{H_{ij \ell}} = x_i -
\e^{\ell} x_j$ and $\alpha_{H_{ij \ell}}^{\vee} = y_i - \e^{-
\ell} y_j$. We have $e_{H_i} =d$ and $\e_{H_i, j} = \frac{1}{d}
\sum_{r=0}^{d-1} \e^{prj} s_i^{pr}$, and $e_{H_{ij \ell}} =2$ and
$ \e_{H_{ij \ell}} = \frac{1}{2} (1+ (-1)^j
\sigma_{ij}^{(\ell)})$. The commutation relation (\ref{commute})
becomes
\begin{multline}\label{commute2} [y_a, x_b] = \delta_{ab} +
\sum_{i=1}^n \delta_{ia}\delta_{ib} \left[ \sum_{j=0}^{d-1}
(\kappa_{j+1} - \kappa_j) \sum_{r=0}^{d-1} \e^{prj} s_i^{pr}
\right] \\ + \sum_{1 \le i < j \le n } \sum_{\ell =0}^{m-1}
\frac{1}{2} (\delta_{ia} - \e^{\ell} \delta_{ja} )(\delta_{ib} -
\e^{-\ell} \delta_{jb}) 2 \kappa_{00} \sigma_{ij}^{\ell}
\end{multline} where $k_{H_{ij \ell}} =  \kappa_{00}$ for all $i,j
, \ell$ and $k_{H_i, j} =  \kappa_j$. We will denote the rational
Cherednik algebra with these parameters by $H_\kappa$.


\subsection{}\label{pbw}
It was proved in \cite[Theorem 1.3]{EtGi} that $H_\kappa$
satisfies a PBW--property, that is, it is isomorphic as a vector
space to $\C[\h] \otimes \C W \otimes \C[\h^*]$ via the
multiplication map.

\subsection{}
Let $H_\kappa (G(m,p,n))$ denote the rational Cherednik algebra of
$G(m,p,n)$ with parameters $\kappa$. We will make considerable use
of the fact that there is an embedding $$H_\kappa (G(m,p,n))
\hookrightarrow H_\mu (G(m,1,n))$$ for an appropriate choice of
$\mu$. This observation is essentially due to Dunkl and Opdam
\cite{DO}.
\begin{thrm}\label{inclusion}
Given $\kappa = (\kappa_{00}, \kappa_1, \ldots, \kappa_{d-1})$,
define $\mu= (\mu_{00}, \mu_1, \ldots, \mu_{d-1}, \mu_d , \ldots,
\mu_{m-1})$ by $\mu_{00}= \kappa_{00}$, $\mu_0=0$, $\mu_i =
\kappa_i/p$, $1 \le i \le d-1$, and $\mu_{sd+t} = \mu_t$ for $1
\le s \le p-1$ and $1 \le t \le d-1$. Then $H_\kappa(G(m,p,n))$ is
the subalgebra of $H_\mu (G(m,1,n))$ generated by $\h$, $\h^*$,
$\sigma_{ij}^{(\ell)}$ for all $i,j,\ell$, and $s_i^p$, $1 \le i
\le n$.
\end{thrm}
\begin{proof}
Write $H_\kappa:= H_\kappa(G(m,p,n))$ and $H_\mu :=
H_\mu(G(m,1,n))$. We need to check that the copies of $\h$, $\h^*$
in $H_\mu$ obey the commutation relations for $H_\kappa$. It is
simply a question of substituting the $\mu$ values into
(\ref{commute2}) above, with $p=1$. We obtain, in $H_\mu$,
\begin{multline*} [y_a, x_b] = \delta_{ab} + \sum_{i=1}^n
\delta_{ia}\delta_{ib} \left[ \sum_{j=0}^{m-1} (\mu_{j+1} - \mu_j)
\sum_{r=0}^{m-1}
\e^{rj} s_i^{r} \right] \\
+ \sum_{1 \le i < j \le n } \sum_{\ell =0}^{m-1} (\delta_{ia} -
\e^{\ell} \delta_{ja} )(\delta_{ib} - \e^{-\ell} \delta_{jb})
\mu_{00} \sigma_{ij}^{\ell} \end{multline*} where $\e:= e^\frac{2
\pi i}{m}$. This may be rewritten as
\begin{multline*} [y_a, x_b] = \delta_{ab} + \sum_{1 \le i < j
 \le n } \sum_{\ell =0}^{m-1} (\delta_{ia} -
\e^{\ell} \delta_{ja} )(\delta_{ib} - \e^{-\ell} \delta_{jb})
\kappa_{00} \sigma_{ij}^{\ell}  + \sum_{i=1}^n
\delta_{ia}\delta_{ib} \frac{1}{p} \sum_{q=0}^{m-1} x_q s_i^q
\end{multline*}
where
$$x_q = \sum_{j=0}^{d-1} \e^{qj} (\kappa_{j+1}- \kappa_j) +
\sum_{j=d}^{2d-1} \e^{qj} (\kappa_{[j+ 1]}-\kappa_{[j]}) + \cdots
+ \sum_{j= (p-1)d}^{pd-1} \e^{qj} (\kappa_{[j+1]} - \kappa_{[j]}
).$$ where $[j]$ denotes the remainder modulo $d$. Write $q =
ap+b$, $0 \le a \le d-1$, $0 \le b \le p-1$. Then
\begin{align*}
x_q &= \sum_{j=0}^{d-1} (\e^{qj}+ \e^{q(j+d)} + \cdots +
\e^{q(j+(p-1)d)}) (\kappa_{j+1} - \kappa_j ) \\
&= \sum_{j=0}^{d-1} \sum_{r=0}^{p-1} \e^{q(j+rd)}( \kappa_{j+1} -
\kappa_j ) \\ &= \sum_{j=0}^{d-1} \alpha_j (\kappa_{j+1}-\kappa_j)
\end{align*}
where $\alpha_j = \e^{qj} \sum_{r=0}^{p-1} \e^{qrd} = \e^{qj}
\sum_{r=0}^{p-1} (e^\frac{2 \pi i}{p} )^{br}$. So
\begin{equation*}\alpha_j = \begin{cases}
\e^{apj} & \text{if $q=ap$,} \\
 0 & \text{if not.} \end{cases}
\end{equation*}
And so
\begin{equation*}
x_q = \sum_{j=0}^{d-1} \alpha_j (\kappa_{j+1}-\kappa_j) =
\begin{cases} 0 & \text{if $q \neq ap$,} \\ p \sum_{j=0}^{d-1}
\e^{apj} (\kappa_{j+1} - \kappa_j) & \text{$q= ap$.} \end{cases}
\end{equation*}
So in $H_\mu$ we have
\begin{multline*} [y_a, x_b] = \delta_{ab} + \sum_{i=1}^n
\delta_{ia}\delta_{ib}\sum_{\stackrel{q=0}{q=ap}}^{m-1} \sum_{j=0}^{d-1} \e^{apj}
(\kappa_{j+1} - \kappa_j) s_i^q \\
+ \sum_{1 \le i < j \le n } \sum_{\ell =0}^{m-1} (\delta_{ia} -
\e^{\ell} \delta_{ja} )(\delta_{ib} - \e^{-\ell} \delta_{jb})
\kappa_{00} \sigma_{ij}^{\ell}. \end{multline*} The first term on
the right hand side is $\sum_{i=1}^n \delta_{ia}\delta_{ib}
\sum_{r=0}^{d-1}\sum_{j=0}^{d-1} \e^{prj} (\kappa_{j+1}- \kappa_j)
s_i^{rp}$. So $y_a, x_b$ obey the commutation relation for
$H_\kappa$.

To finish the proof of the theorem, we may define a map
$$T(\h \oplus \h^*) \ast G(m,p,n) \rightarrow H_\mu$$
in the obvious way. We have checked above that the commutation
relations for $H_\kappa$ are in the kernel, and it is easily seen
that the other relations for $H_\kappa$ are in the kernel as well.
Thus, there is a well-defined map
$$H_\kappa \rightarrow H_\mu.$$
To check that this is injective, consider an element of $H_\kappa$
which is mapped to zero. Write it in terms of a PBW-basis of
$H_\kappa$, and observe that all the coefficients must therefore
be 0, since a PBW-basis of $H_\kappa$ is mapped into a subset of a
PBW-basis of $H_\mu$.
\end{proof}

\subsection{The Dunkl representation}\label{dunkl}
It is well--known (see for instance, \cite{DO}, \cite[Proposition
4.5]{EtGi}) that $H_\kappa$ acts on $\C [\h] \otimes
\textsf{triv}$ where $\textsf{triv}$ denotes the trivial
representation of $W$. Furthermore, this action is faithful, and
if $\C [\h] \otimes \textsf{triv}$ is identified with $\C [\h]$,
then the action of $y \in \h$ is given by a
differential--difference operator called a Dunkl operator:
$$T_y = \partial_y + \sum_{H \in \mathcal{A}}
 \frac{\langle \alpha_H , y \rangle }{\alpha_H} \sum_{i=1}^{e_H -1} e_H k_{H, i} \e_{H, i}$$
If $\h^{\mathrm{reg}} = \h \setminus \cup_{H \in \mathcal{A}} H$
then the Dunkl representation defines an injective homomorphism
$$H_\kappa \hookrightarrow \mathcal{D}(\h^{\mathrm{reg}})*W$$
called the Dunkl representation. If $\delta = \prod_{H \in
\mathcal{A}} \alpha_H \in \C[\h]$, then $\C [\h^{\mathrm{reg}}] =
\C [\h]_\delta$ and the induced map
$$H_\kappa|_{\h^{\mathrm{reg}}}:=H_\kappa \otimes_{\C[\h]} \C [\h^{\mathrm{reg}}]
\rightarrow \mathcal{D}(\h^{\mathrm{reg}})*W$$
is an isomorphism (\cite[Theorem 5.6]{GGOR}).

\subsection{} Let $\theta_\mu : H_\mu \rightarrow
\mathcal{D}(\h^{\mathrm{reg}}) \ast G(m,1,n)$ be the Dunkl
representation of $H_\mu$ and let $\theta_\kappa: H_\kappa
\rightarrow \mathcal{D}(\h^{\mathrm{reg}}) \ast G(m,p,n)$ be the
Dunkl representation of $H_\kappa$. Regarding $H_\kappa$ as a
subalgebra of $H_\mu$, we wish to show that $\theta_\mu
|_{H_\kappa} = \theta_\kappa$. For this, it suffices to check that
$\theta_\mu (y) = \theta_\kappa(y)$ for all $y \in \h$. But
$\theta_\mu(y)$ and $\theta_\kappa (y)$ may be regarded as
differential-difference operators acting on the polynomial ring
$\C[\h]$, so it suffices to check that their values on polynomials
are the same. If $p \in \C[\h]$ then $\theta_\mu(y) (p) = y \cdot
p \otimes 1$ where we identify $p \in \C[\h]$ with $p \otimes 1
\in \C[\h] \otimes \mathsf{triv}$. But $y \cdot p \otimes 1 =
[y,p]_\mu \otimes 1$ where $[y,p]_\mu$ denotes the commutator in
$H_\mu$. This may be written in terms of commutators $[y,x]_\mu$
for $x \in \h^*$. But $[y,x]_\mu = [y,x]_\kappa$ for all $y \in
\h$, $x \in \h^*$, where $[y,x]_\kappa$ denotes the commutator in
$H_\kappa$. So $[y,p]_\mu \otimes 1 = [y,p]_\kappa \otimes 1 =
\theta_\kappa(y) (p)$. So $\theta_\kappa(y) = \theta_\mu (y)$ as
required. We have proved the following lemma.
\begin{lem}\label{dunklrep}
If $H_\kappa$ is the rational Cherednik algebra of $G(m,p,n)$, and
we consider $H_\kappa$ as a subalgebra of $H_\mu$ as in Theorem
\ref{inclusion}, then the  Dunkl representation $\C[\h] \otimes
\triv$ of $H_\mu$ restricts to the Dunkl representation of
$H_\kappa$. \end{lem}

\section{Category $\mathcal{O}$}\label{cato}
\subsection{}
In this section, we will review the theory for a general complex
reflection group $W$ and its rational Cherednik algebra $H_\kappa$
depending on some collection of complex parameters $\kappa =
(k_{H,i})_{H \in \mathcal{A}, \: 0 \le i \le e_H-1}$.

\subsection{}
Following \cite{BEG1}, let $\mathcal{O}$ be the abelian category of
finitely-generated $H_\kappa$--modules $M$ such that
for $P \in \C[\h^*]^W$, the action of $P - P(0)$ is locally
nilpotent. Let $\textsf{Irrep}(W)$ denote the set of isoclasses of
simple $W$--modules. Given $\tau \in \irrep(W)$, define the
standard module $M(\tau)$ by:
$$M(\tau) = H_{\kappa} \otimes_{\C [\h^*]*W} \tau$$
where for $p \in \C[\h^*], w \in W$ and $v \in \tau$, $pw \cdot v
:= p(0) wv$.

\subsection{}
In \cite{DO}, it is proved that $M(\tau)$ has a unique simple
quotient $L(\tau)$, and \cite{GGOR} prove that $\{ L(\tau) | \tau
\in \irrep(W) \}$ is a complete set of nonisomorphic simple
objects of $\cO$, and that every object of $\cO$ has finite
length.

\subsection{}
\label{zdef} By \cite{GGOR}, if $z := \sum_{H \in \mathcal{A}}
\sum_{i=1}^{e_H -1} e_H k_{H,i} \e_{H, i} \in Z(\C W)$, and
$\mathfrak{d} := \sum_i x_i y_i \in H_\kappa$, then
$\textsf{eu}_\kappa = \textsf{eu}:= \mathfrak{d} -z$ has the
property that $[\textsf{eu},x]=x$ for all $x \in \h^*$ and
$[\textsf{eu},y] = -y$ for all $y \in \h$ and $[\textsf{eu},w]=0$
for all $w \in W$. The action of $\textsf{eu}$ on $M(\tau)$ is
diagonalisable and the eigenspaces are $\C [\h]_d \otimes \tau$,
$d \ge 0$, where $\C[\h]_d$ denotes the homogeneous polynomials in
$\C [\h]$ of degree $d$. The eigenvalue of $\mathsf{eu}$ on
$\C[\h]_d \otimes \tau$ is $d- \theta(z)$ where $\theta(z)$ is the
eigenvalue of $z$ on $\tau$. In particular, the lowest eigenvalue
of $\mathsf{eu}$ on $M(\tau)$ is $-\theta(z)$.

\subsection{}\label{catoaltdef}
A useful alternative definition of category $\cO$ is quoted in
\cite[Section 2.1]{CE}. Category $\cO$ may be defined as the
category of $H_\kappa$--modules $V$ such that $V$ is a direct sum
of generalised $\mathsf{eu}_\kappa$ eigenspaces, and such that the
real part of the spectrum of $\mathsf{eu}_\kappa$ is bounded
below. It is clear from this definition that every
finite-dimensional $H_\kappa$--module belongs to $\cO$.

\subsection{}
The group $B_W:= \pi_1 (\h^{\mathrm{reg}}/W)$ is called the braid
group of $W$. In \cite{GGOR}, a functor
$$\texttt{KZ} : \cO \rightarrow \C B_W-\mathrm{mod}$$
is constructed as follows: If $M \in \cO$ then
$M|_{\h^{\mathrm{reg}}}:= \C [\h^{\mathrm{reg}}] \otimes_{\C[\h]}
M$ is a finitely-generated module over $\C [\h^{\mathrm{reg}}]
\otimes_{\C[\h]} H_\kappa \cong \mathcal{D}(\h^{\mathrm{reg}})*W$.
In particular, $M$ is a $W$--equivariant $\mathcal{D}$--module on
$\h^{\mathrm{reg}}$ and hence corresponds to a $W$--equivariant
vector bundle on $\h^{\mathrm{reg}}$ with a flat connection
$\nabla$. The horizontal sections of $\nabla$ define a system of
differential equations on $\h^{\mathrm{reg}}$ which, by a process
described in \cite{BMR} and \cite{rouquiersurvey}, give a
monodromy representation of $\pi_1 (\h^{\mathrm{reg}}/W)$. By
definition, $\texttt{KZ}(M)$ is the monodromy representation of
$\pi_1 (\h^{\mathrm{reg}}/W)$ associated to $M$.

\subsection{}
\label{kznumbers} By \cite[4.12]{BMR} and \cite[Section
5.25]{GGOR}, the monodromy representation factors through the
Hecke algebra $\mathcal{H}$ of $W$. This is the quotient of $\C
B_W$ by the relations:
$$(T-1) \prod_{j=1}^{e_H-1} (T- \det(s)^{-j} e^{-2 \pi i k_{H,j}})$$
for $H \in \mathcal{A}$, $s \in W$ the reflection around $H$ with
nontrivial eigenvalue $e^{2 \pi i /e_H}$, and $T$ an
$s$--generator of the monodromy around $H$. The parameters differ
from those given in \cite{GGOR} because the idempotent $\e_j(H)$
of \cite{BMR} is the $\e_{-j , H}$ of \cite{GGOR}.

\subsection{}\label{otor}
Therefore, $\texttt{KZ}$ gives a functor $\texttt{KZ}: \cO
\rightarrow \mathcal{H}-\mathrm{mod}$. By \cite[Section
5.3]{GGOR}, $\texttt{KZ}$ is exact, and if $\cO_{\mathrm{tor}}$ is
the full subcategory of those $M$ in $\cO$ such that
$M|_{\h^{\mathrm{reg}}} =0$ then $\texttt{KZ}$ gives an
equivalence $\cO/\cO_{\mathrm{tor}} \tilde{\rightarrow}
\mathcal{H}-\mathrm{mod}$ \cite[Theorem 5.14]{GGOR}.

\section{The Hecke algebra of $G(m,p,n)$}\label{hecke}
We now identify the algebra $\mathcal{H} =
\mathcal{H}(\kappa_{00}, \kappa_1, \ldots, \kappa_{d-1})$ through
which $\texttt{KZ}$ factors, in the case of $W= G(m,p,n)$. By
\cite[Prop 4.22]{BMR}, this algebra is generated by $(T_s)_{s \in
\mathcal{N}(D)}$ where $\mathcal{N}(D)$ is the set of nodes of the
braid diagram $D$ of $W$. The generators $T_s$ are subject to the
braid relations defined by $D$, together with the relations of
\ref{kznumbers}. From the braid diagram in \cite[Table 1]{BMR} we
see that if $p>2$ then $\mathcal{H}$ is generated by $T_s,
T_{t_2}, T_{t^{'}_2}, T_{t_3} , \ldots, T_{t_n}$ subject to
 the following relations
\begin{align*}
T_s T_{t^{'}_2} T_{t_2} - T_{t^{'}_2} T_{t_2} T_s &=0 \\
T_{t^{'}_2} T_{t_3} T_{t^{'}_2} - T_{t_3} T_{t^{'}_2} T_{t_3} &=0
\\
T_{t_2} T_{t_3} T_{t_2} - T_{t_3} T_{t_2} T_{t_3} &=0 \\
T_{t_3} T_{t^{'}_2} T_{t_2} T_{t_3} T_{t^{'}_2} T_{t_2} -
T_{t^{'}_2} T_{t_2} T_{t_3} T_{t^{'}_2} T_{t_2} T_{t_3} &=0 \\
T_{t_2} T_s (T_{t^{'}_2} T_{t_2})_{p-1} - T_s (T_{t^{'}_2}
T_{t_2})_p &=0 \\
[T_{t_i}, T_s] &=0 & i \ge 3\\
[T_{t_2}, T_{t_i}] &=0 & i \ge 4\\
[T_{t^{'}_2}, T_{t_i}] &=0 & i \ge 4\\
[T_{t_i}, T_{t_j}] &=0 & i, j \ge 3, &\: |i-j| \ge 2\\
T_{t_i} T_{t_{i+1}} T_{t_i} - T_{t_{i+1}} T_{t_i} T_{t_{i+1}} &=0
& i \ge 3 \\
(T_s -1) \prod_{j=1}^{d-1} (T_s - \e^{-pj} e^{2 \pi i \kappa_j} )
&= 0 \\
(T_{t_i} - 1)(T_{t_i} + e^{2 \pi i \kappa_{00}} ) &=0 & \forall i
\\ (T_{t_2^{'}} -1)(T_{t_2^{'}} + e^{2 \pi i \kappa_{00}}) &= 0 \\
\end{align*}
where if $x, y$ are generators then $(xy)_r$ denotes the word
$(xy)^{r/2}$ if $r$ is even or $(xy)^{(r-1)/2}x$ if $r$ is odd. If
$p=2$ then we see from \cite[Table 2]{BMR} that $\mathcal{H}$ is
the algebra described above, except that the relation $T_{t_3}
T_{t^{'}_2} T_{t_2} T_{t_3} T_{t^{'}_2} T_{t_2} - T_{t^{'}_2}
T_{t_2} T_{t_3} T_{t^{'}_2} T_{t_2} T_{t_3} =0$ is omitted.

We wish to verify that the algebra presented by these generators
and relations is the same as the Hecke algebra of $G(m,p,n)$ as
defined in \cite[2.A]{GJ}. We set $a_0 = T_s$, $a_1=
-T_{t^{'}_2}$, $a_2 = -T_{t_2}$ and $a_i = -T_{t_i}$ for $i \ge
3$. In the $p=2$ case, we see that we have exactly the relations
of \cite[2.A]{GJ}, except that \cite{GJ} have the additional
relation $(a_1 a_2 a_3)^2 = (a_3 a_1 a_2)^2$. However, this
additional relation follows from the other relations, since $(T_s
-1) \prod_{j=1}^{d-1} (T_s - \e^{-pj} e^{2 \pi i \kappa_j} ) = 0$
implies that $a_0 = T_s$ is invertible, and we can then check that
$(a_1 a_2 a_3)^2 a_0 = (a_3 a_1 a_2)^2 a_0$, using the relations
listed above.

If $p>2$, we  again have the same relations as \cite[2.A]{GJ},
except that the relation
\begin{equation}\label{rel1}a_0 a_1 a_2 = (q^{-1} a_1 a_2)^{2-p} a_2
a_0 a_1 + (q-1) \sum_{k=1}^{p-2} (q^{-1} a_1 a_2)^{1-k} a_0
a_1\end{equation} of \cite{GJ} has been replaced by the relation
\begin{equation}\label{rel2}a_2 a_0 (a_1 a_2)_{p-1} = a_0 (a_1
a_2)_p.\end{equation} An explicit calculation shows that in the
presence of the other relations, (\ref{rel1}) and (\ref{rel2}) are
equivalent, where $q= e^{2 \pi i \kappa_{00}}$. We therefore
obtain the following lemma.
\begin{lem}
The Hecke algebra $\mathcal{H}$ through which the functor
$\texttt{KZ}$ factors is isomorphic to the Hecke algebra denoted
$\mathfrak{H}^{q,x}_{m,p,n} (\mathbb{C})$ in \cite{GJ}, with
parameters $q= e^{2 \pi i \kappa_{00}}$ and $x_1 = 1$, $x_j =
\e^{-p(j-1)} e^{-2 \pi i \kappa_{j-1}}$ for $j >1$.
\end{lem}

We will need some facts about the representation theory of this
algebra. Specifically, we will need to use the parametrisation of
the simple modules that is the main result of \cite[Section
3]{GJ}.

\subsection{}
As in \cite[2.B]{GJ}, we make the following definitions. For $1\le
i \le m$, write $i=sp+t$, with $0 \le s \le d-1$ and $1 \le t \le
p$. Let $\eta_p = e^{\frac{2 \pi i}{p}} = \e^d$ and let $Q_i =
\eta_p^{t-1}y_{s+1}$ where $y_{s+1}$ is chosen so that $y_{s+1}^p
= x_{s+1}$. In this way, we get a new sequence of complex numbers
$Q:= (Q_1, \ldots , Q_m)$. Now we follow \cite[2.C]{GJ}. Let
$\Pi^m_n$ denote the set of multipartitions $\lambda =
(\lambda^{(i)})_{1 \le i \le m}$ of $n$ with $m$ parts. Then there
is a permutation $\varpi$ that acts on $\Pi^m_n$ as follows.
The permutation $\varpi$ may be expressed in
cycle notation as
$$\varpi = (1,2, \ldots, p)(p+1, p+2, \ldots 2p) \cdots ((d-1)p
+1, \ldots dp).$$ The action of $\varpi$ on a multipartition
$\lambda = (\lambda^{(i)})$ is defined by $\varpi (\lambda)^{(i)}
= \lambda^{\varpi^{-1} (i)}$. Let $\mathcal{L}$ be a set of
representatives of the orbits of this action of $\varpi$ on
$\Pi^m_n$ and for $\lambda \in \Pi^m_n$, let $o_\lambda =
\mathrm{min} \{ k \in \mathbb{N}_{>0} | \varpi^k (\lambda) =
\lambda \}$.

\subsection{}
In Theorem \ref{hsimples} below, the set of Kleshchev
multipartitions in $\Pi^m_n$ is defined with respect to the
parameters $q$ and $Q= (Q_1, \ldots, Q_m)$.  To define the set of
Kleshchev multipartitions, we first need the definition of the
\emph{residue} of a node in a multipartition $(\lambda^{(1)},
\ldots \lambda^{(m)} )$. If $x$ is a node in column $j(x)$ and row
$i(x)$ of $\lambda^{(k)}$, we define the residue $\res(x) = Q_k
q^{j(x) -i(x)}$. We say $y \notin \lambda$ is an \emph{addable}
$a$--node if $\lambda \cup \{ y\}$ is a multipartition and
$\res(y)= a$. We say $y \in \lambda$ is a \emph{removable}
$a$--node if $\lambda \setminus \{ y\}$ is a multipartition and
$\res(y)=a$. A node $x \in \lambda^{(i)}$ is said to be
\emph{below} a node $y \in \lambda^{(j)}$ if either $i>j$ or else
$i=j$ and $x$ is in a lower row than $y$. A removable node $x$ is
called a \emph{normal} $a$--mode if whenever $y$ is an addable
$a$--node which is below $x$ then there are more removable
$a$--nodes between $x$ and $y$ than there are addable $a$--nodes.
A removable $a$--node is called \emph{good} if it is the highest
normal $a$--node of $\lambda$. The set of Kleshchev
multipartitions is defined inductively by declaring that the empty
multipartition is Kleshchev, and that a multipartition $\lambda$
is Kleshchev if and only if there is a node $y \in \lambda$ which
is a good $a$--node, for some $a$, such that $\lambda \setminus \{
y\}$ is Kleshchev. A more detailed exposition, including examples,
may be found in the introduction to the paper \cite{arikimathas}.

By the remark following \cite[Theorem 3.1]{GJ}, we have the
following theorem.
\begin{thrm}\cite{GJ}\label{hsimples}
Suppose $q$ is not a root of unity. Then the set of simple
$\mathfrak{H}^{q,x}_{m,p,n}$--modules is in bijection with the set
$\{ (\lambda, i) | \lambda \in \Lambda^{0} \cap \mathcal{L} \:
\mathrm{and} \: i \in [0, \frac{p}{o_\lambda} -1]\}$ where
$\Lambda^0$ denotes the set of \emph{Kleshchev} multipartitions in
$\Pi^m_n$, and $\mathcal{L}$ is defined as above.
\end{thrm}

\subsection{}
In this paper, we will be primarily interested in values of the
parameters such that the equation
\begin{equation}\label{kapparels} d \kappa_1 + m (n-1)\kappa_{00}
= -1-m(n-1)-d \end{equation} holds. If this is the case, then
$$q^{p(n-1)} = \e^{-p} e^{-2 \pi i \kappa_1} = x_2.$$ We will use
the parametrisation given by Theorem \ref{hsimples} to work out
how many simple modules $\mathcal{H}$ has for generic choices of
$\kappa_{00}, \kappa_1, \ldots , \kappa_{d-1}$ satisfying
(\ref{kapparels}). This will be used in the next section to give
an upper bound on the number of finite-dimensional simple objects
in category $\cO$.
\section{Simple modules for the Hecke algebra}
\subsection{}
In this section, we have the standing assumptions that the
parameters are $q= e^{2 \pi i \kappa_{00}}$ and $Q= (Q_1, \ldots,
Q_m)$ where the $Q_i$ are defined as above. In particular,
$Q_{sp+t} = \eta_p^{t-1} x_{s+1}^{\frac{1}{p}}$. But $x_{s+1} =
\e^{-ps} e^{-2 \pi i \kappa_{s}}$ so $Q_{sp+t} = \eta_p^{t-1}
\e^{-s} e^{-2 \pi i \kappa_s /p} = \e^{d(t-1) -s} e^{-2 \pi i
\kappa_s /p}$. In particular, $Q_1 = 1$, $Q_2 = \eta_p$, $\ldots$,
$Q_p = \eta_p^{p-1}$ and $Q_{p+1} = y_2 = q^{n-1}$, $Q_{p+2} =
\eta_p q^{n-1}$, $\ldots$, $Q_{2p} = \eta_p^{p-1} q^{n-1}$. We
claim that, for these values of $Q_i$ and $q$, and for a generic
choice of the $\kappa_{00}$ and $\kappa_i$, there are exactly $p$
multipartitions in $\Pi^m_n$ which are not Kleshchev and they can
be described as follows.

Let $\rho = \begin{matrix} {\tiny\yng(3)} & \cdots & \tiny\yng(1)
\end{matrix}$, a partition of $n$, and for $1 \le i \le m$, define $\rho_i
\in \Pi^m_n$ to be the multipartition with $\rho$ in the
$i^{\mathrm{th}}$ place and $\O$ everywhere else. Then we have the
following lemma.
\begin{lem}\label{kleshchev}
With the above choices of $Q_i$ and $q$, the non-Kleshchev
multipartitions in $\Pi^m_n$ are prescisely the $\rho_i$ where $1
\le i \le p$.
\end{lem}
\begin{proof}
First, note that for a generic choice of the parameters,
$\kappa_{00} \notin \mathbb{Q}$ and so we may assume that $q$ is
not a root of unity.

 Let $\lambda = (\lambda^{(1)}, \ldots,
\lambda^{(m)}) \in \Pi^m_n$. Suppose that $\lambda \neq \rho_i$
for $1 \le i \le p$. We must show that $\lambda$ is Kleshchev. We
will show that we can repeatedly remove good nodes from $\lambda$
until we reach the empty partition. First, let $i
>2p$. We will show that we can reduce to the case $\lambda^{(i)}
=\O$. Recall that $Q_i$ is of the form $y_a \eta_p^b$ for some $a$
and $b$. Thus, the residue of a node in $\lambda^{(i)}$ is of the
form $q^c y_a \eta_p^b$ for some $a,b,c$. If this is equal to the
residue of a node in some $\lambda^{(j)}$ where $j \neq i$, then
we must have $q^c y_a \eta_p^b = q^{c^{'}} y_{a^{'}}
\eta_p^{b^{'}}$ for some $a^{'}, b^{'}, c^{'}$. So $q^{c- c^{'}} =
y_{a^{'}}y_a^{-1} \eta_p^{b^{'}-b}$. Hence, $e^{2 \pi i
\kappa_{00} (c-c^{'})} = \e^{-(a^{'}-1)} e^{-2 \pi i
\kappa_{a^{'}} /p} \e^{(a-1)} e^{2 \pi i \kappa_a /p}
\e^{d(b^{'}-b)}$. So
$$\mathrm{exp} (2 \pi i \kappa_{00} (c-c^{'}) + 2 \pi i \frac{a^{'}-1}{m}
+2 \pi i \frac{\kappa_{a^{'}}}{p} - 2 \pi i \frac{a-1}{m} -2 \pi i
\frac{\kappa_a}{p} - 2 \pi i \frac{d(b^{'}-b)}{m} ) =1.$$ For a
generic choice of the $\kappa_i$ and $\kappa_{00}$, this can only
happen if $c=c^{'}$ and $a=a^{'}$. But then this forces $i=j$.
Hence, the only nodes in $\lambda$ which can have the same residue
as a node in $\lambda^{(i)}$ are the other nodes in
$\lambda^{(i)}$. Let $x$ be the rightmost node in the bottom row
of $\lambda^{(i)}$. Any node with the same residue as $x$ has
residue $q^{j(x)-i(x)} Q_i$ and so must lie on the same diagonal
as $x$. But, by the choice of $x$, there can be no addable or
removable nodes on the same diagonal as $x$. So $x$ is a good
node, and we may remove $x$. Continuing inductively, the nodes of
$\lambda^{(i)}$ may be removed one at a time, and we conclude that
the multipartition with $\lambda^{(i)}$ replaced by $\O$ is
Kleshchev. Note that this is not necessarily a multipartition of
$n$; it is a multipartition of the integer $\sum_{j \neq i}
|\lambda^{(j)} | \le n$.

We are now reduced to checking that all the multipartitions of the
form $$\lambda = (\lambda^{(1)} , \lambda^{(2)}, \ldots ,
\lambda^{(p)} , \ldots ,\lambda^{(2p)}, \O, \O , \ldots, \O)$$
with $\sum_i |\lambda^{(i)}| \le n$ and $\lambda \neq \rho_i$, $1
\le i \le p$, are Kleshchev. We now show that we can reduce to the
case $\lambda^{(2p)} = \lambda^{(2p-1)} = \cdots = \lambda^{(p+1)}
= \O$. First, consider $\lambda^{(2p)}$. Suppose $\lambda^{(2p)}$
has $b$ rows and that the lowest row has length $a$. Let $x$ be
the rightmost node of the bottom row of $\lambda^{(2p)}$. Then the
residue of $x$ is $q^{a-b} \cdot q^{n-1} \eta_p^{p-1} =: r$. Then
$x$ is a normal $r$--node because $\lambda^{(2p)}$ has no addable
$r$--nodes which are below $x$, since such a node would have to
lie on the same diagonal as $x$, which is impossible by choice of
$x$. And if $i>2p$ then $\lambda^{(i)} =\O$ has no addable
$r$--node, since the only node that can be added to
$\lambda^{(i)}$ has residue $Q_i=Q_{sp+t}$ for some $s \ge 2$ and
some $t$. This equals $\e^{d(t-1)-s}e^{-2 \pi i \kappa_s /p}$
which cannot equal $r$ for a generic choice of $\kappa_s$. Hence,
there are no addable $r$--nodes of $\lambda$ which are below $x$,
so $x$ is a normal $r$--node. We must show that $\lambda$ contains
no higher normal $r$--nodes. Certainly $\lambda^{(2p)}$ contains
no higher normal $r$--node, since any $r$--node must lie on the
same diagonal as $x$ and so cannot be removable. Any normal
$r$--node not in $\lambda^{(2p)}$ must lie in $\lambda^{(p)}$,
because no power of $q$ can be equal to a power of $\eta_p$ since
$q$ is not a root of unity. Suppose then that the node $y \in
\lambda^{(p)}$ has residue $r$. Say $y$ lies in row $i(y)$ and
column $j(y)$ of $\lambda^{(p)}$. Then $q^{j(y)-i(y)} \eta_p^{p-1}
= q^{a-b+n-1} \eta_p^{p-1}$. Hence, $j(y)-i(y) = a-b+n-1$. Suppose
$\lambda^{(p)}$ has $d$ columns. Then $a-b+n-1 =j(y)-i(y) \le
d-1$. So $n+1 \le n+a \le b+d$. But $b+d \le |\lambda^{(2p)}|
+|\lambda^{(p)}| \le n$, a contradiction. Hence, no such $y$
exists, and $x$ is the highest normal $r$--node, and so is good.
Removing $x$ and continuing inductively, we may take
$\lambda^{(2p)} = \O$. The same argument shows that we may take
$\lambda^{(p+i)} = \O$, $1 \le i \le p$, as claimed.

We are now reduced to showing that those $\lambda$ of the form
$\lambda = (\lambda^{(1)} , \ldots, \lambda^{(p)} , \O , \ldots
,\O)$ with $\lambda \neq \rho_i$, $1 \le i \le p$, are Kleshchev.
First, consider $\lambda^{(p)}$. Let $b$ be the number of columns
of $\lambda^{(p)}$. Then $b<n$ since $\lambda \neq \rho_p$. If $x
\in \lambda^{(p)}$ then $\mathrm{res}(x) = q^{j(x)-i(x)}
\eta_p^{p-1} $ where $j(x)-i(x) \le b -1$. Let $x$ be the
rightmost node in the bottom row of $\lambda^{(p)}$ and let $r=
\mathrm{res}(x)$. Then $x$ is the only removable $r$--node in
$\lambda$ since there can be no $r$--nodes in $\lambda^{(i)}$ with
$i \neq p$, again because $q$ is not a root of unity. Furthermore,
the only way there can be an addable $r$--node below $x$ is if
such a node can be added to the empty diagram $\lambda^{(2p)}$.
Such a node would have residue $q^{n-1} \eta_p^{p-1}$ and so we
would have to have $j(x)-i(x) = n-1 \le b-1$. This contradicts
$b<n$. Hence, $x$ is a good node and may be removed. Continuing
inductively, we may remove all the nodes in $\lambda^{(p)}$. The
same argument works for all the $\lambda^{(i)}$, $1 \le i \le p$,
and so we can get to the empty partition from $\lambda$ by
successively removing good nodes. Hence, $\lambda$ is Kleshchev.

To complete the proof, we observe that $\rho_i$ is not Kleshchev
for $1 \le i \le p$, since the only removable node in $\rho_i$ is
the node at the end of the row $\rho_i^{(i)}$. This node has
residue $q^{n-1} \eta_p^{i-1}$. It is not normal as there is an
addable node in $\rho_i$ which is the unique node that may be
added to the empty diagram $\rho_i^{(p+i)}$, and this node also
has residue $q^{n-1} \eta_p^{i-1}$.
\end{proof}

\begin{cor}\label{numbersimples}
For generic values of the parameters satisfying $d \kappa_1 +
m(n-1) \kappa_{00} = -1-m(n-1)-d$, the algebra $\mathcal{H}$ has
$|\irrep(W)|-1$ nonisomorphic simple modules.
\end{cor}
\begin{proof}
By Theorem \ref{hsimples}, the simples are indexed by pairs $
(\lambda,i)$ such that $\lambda \in \Lambda^0 \cap \mathcal{L}$
and $0 \le i \le \frac{p}{o_\lambda}-1$. Let $T$ be the set of all
such pairs. Note that $\{ \rho_1, \rho_2, \ldots, \rho_p\}$ is an
orbit of $\varpi$ on $\Pi^m_n$. Let $S$ be the set of all pairs
$(\lambda,i)$ with $\lambda \in \mathcal{L}$ and $0 \le i \le
\frac{p}{o_\lambda} -1$. Choose $\rho_1$ to be the representative
of the orbit of $\rho_1$ in $\mathcal{L}$. Then $(\rho_1,0) \in S$
and $(\rho_1,0) \notin T$. Furthermore, if $\lambda \in
\mathcal{L}$ and $\lambda \neq \rho_1$ then $\lambda$ is Kleshchev
by \ref{kleshchev} and so $(\lambda,i) \in T$ for all $0 \le i \le
\frac{p}{o_\lambda}-1$. Hence, $T = S \setminus \{(\rho_1, 0)\}$
and so $|T| = |S|-1$. It remains to show that $|S|= |\irrep(W)|$.

By \cite[Theorem 2.6]{arikigrpn}, $|S|$ is the number of simple
modules for the generic version $K \mathcal{H}$ of the Hecke
algebra $\mathcal{H}$, defined over the field $K=\mathbb{C}(q,x_1,
\ldots, x_d)$ with $q$ and the $x_i$ being indeterminates. But by
the proof of \cite[Theorem 4.24]{BMR}, $K \mathcal{H}$ is
isomorphic to the group algebra $KW$ and so $|S|=|\irrep(W)|$.
\end{proof}

\subsection{Application to category $\cO$}
We may use Corollary \ref{numbersimples}, together with the
functor $\KZ$, to get some information about the category $\cO$ of
$H_\kappa$--modules.
\begin{thrm}\label{infdim}
Let $\kappa_i, \kappa_{00}$ be chosen generically so that $d
\kappa_1 +m(n-1) \kappa_{00} = -1-m(n-1)-d$. Then there is exactly
one finite-dimensional simple module in category $\cO$, namely
$L(\mathsf{triv})$.
\end{thrm}
\begin{proof}
First, recall from Section \ref{otor} that $\KZ$ is an exact
functor, and that every $\mathcal{H}$--module is the image of some
object of $\cO$ under $\KZ$. Suppose $X$ is a simple
$\mathcal{H}$--module. Then, using exactness of $\KZ$, we can find
a simple object $L(\tau) \in \cO$ such that $\KZ(L(\tau)) = X \neq
0$. Since $\KZ(L(\tau)) \neq 0$, we get
$L(\tau)|_{\h^{\mathrm{reg}}} \neq 0$ and hence $\dim (L(\tau)) =
\infty$. Since each of the $|\irrep(W)|-1$ simple
$\mathcal{H}$--modules is then the image of some
infinite-dimensional $L(\tau)$, we conclude that at least
$|\irrep(W)|-1$ of the $L(\tau)$ are infinite-dimensional. (This
argument is based on \cite[Lemma 3.11]{BEG2}).

It remains to show that $L(\mathsf{triv})$ is finite-dimensional.
It was shown in Lemma \ref{dunklrep} that the Dunkl representation
$\C[\h] = M(\mathsf{triv})$ of $H_\kappa$ is the restriction to
$H_\kappa$ of the Dunkl representation of the larger algebra
$H_\mu$, where $H_\mu$ is the rational Cherednik algebra for
$G(m,1,n)$ with parameters $\mu_{00} = \kappa_{00}$ and $\mu_i =
\kappa_i/p$. It follows from the hypothesis on $\kappa$ that
$$ m \mu_1 + m(n-1) \mu_{00}= -1-m(n-1)-d.$$
In the terminology of \cite{CE}, this says that $(\mu_{00}, \mu_1,
\ldots, \mu_{m-1})$ belongs to the set $E_r$ where $r= m(n-1) +d+1
= m(n-1) +q$ where $q=d+1$ and $1 \le q \le m-1$. We are now in a
position to apply \cite[Proposition 4.1]{CE}. Let $\h_q$ be the
representation of $G(m,1,n)$ with $\h_q = \C^n$ as a vector space,
and on which $S_n$ acts by permuting the coordinates, and $s_i$
acts by multiplying the $i^{\mathrm{th}}$ coordinate by
$\e^{-d-1}$. Then \cite[Proposition 4.1]{CE} states that the
polynomial representation $M(\triv)$ of $H_\mu$ contains a copy of
$\h_q$ in degree $m(n-1) +d+1$ consisting of singular vectors
(vectors which are killed by $\h \subset H_\mu$). Furthermore, if
we set $\tilde{Y}_c$ to be the quotient of $M(\mathsf{triv})$ by
the ideal generated by this copy of $\h_q$, then \cite[Theorem
4.3(i)]{CE} states that $\tilde{Y}_c$ is finite-dimensional
provided $\kappa_{00} \notin \mathbb{Q}$. Hence, there is an exact
sequence of $H_\mu$--modules
$$M(\triv) \rightarrow \tilde{Y}_c \rightarrow 0.$$
Since $H_\kappa \subset H_\mu$, these are also $H_\kappa$--module
maps, and so the Dunkl representation $M(\triv)$ of $H_\kappa$ has
a finite-dimensional quotient. It follows that the unique smallest
quotient $L(\triv)$ of the $H_\kappa$--module $M(\triv)$ is
finite-dimensional, as required.
\end{proof}

\subsection{}
This ends our study of the Hecke algebra and category $\cO$. We
now wish to apply Theorem \ref{infdim} to study a special object
in $\cO$ whose associated graded module will yield the desired
quotient of the ring of diagonal coinvariants.

\section{Shifting}\label{shifting}
\subsection{} In this section we will construct a one-dimensional
$H_\kappa$--module $\Lambda$ for particular values of $\kappa$,
and then construct a shifted version of $\Lambda$ which will be a
finite-dimensional $H_\kappa$--module $L$. In the next section we
will show how $L$ is related to the diagonal coinvariants of $W$.

\subsection{A one-dimensional module}
\begin{thrm}\label{onedim}
Suppose we choose the parameters $\kappa_{00}, \kappa_i$ such that
$d \kappa_1 + m(n-1) \kappa_{00} = -1$. Then $H_\kappa$ has a
one-dimensional module $\Lambda$ on which $\h$ and $\h^*$ act by
zero and $W=G(m,p,n)$ acts by the trivial representation.
\end{thrm}
\begin{proof}
Let $\Lambda$ be the trivial $W$--module. We can make $\Lambda$
into a $T(\h \oplus \h^*) \ast W$--module by making $\h$, $\h^*$
act by 0. So we are reduced to showing that the defining relations
of $H_\kappa$ act by $0$ on $\Lambda$. The only relation that may
cause difficulty is the commutation relation of section
\ref{commutationrel}. We need to check that the commutation
relation between $y_a$ and $x_b$ acts by $0$ on $\Lambda$ for all
$a,b$. Recall that
\begin{multline}\label{commute3} [y_a, x_b] = \delta_{ab} +
\sum_{i=1}^n \delta_{ia}\delta_{ib} \left[ \sum_{j=0}^{d-1}
(\kappa_{j+1} - \kappa_j) \sum_{r=0}^{d-1} \e^{prj} s_i^{pr}
\right] \\ + \sum_{1 \le i < j \le n } \sum_{\ell =0}^{m-1}
\frac{1}{2} (\delta_{ia} - \e^{\ell} \delta_{ja} )(\delta_{ib} -
\e^{-\ell} \delta_{jb}) 2 \kappa_{00} \sigma_{ij}^{\ell}
\end{multline}
First, if $a \neq b$ then the right hand side of (\ref{commute3})
becomes
$$\sum_{1 \le i < j \le n } \sum_{\ell =0}^{m-1}
 (\delta_{ia} - \e^{\ell} \delta_{ja} )(\delta_{ib} - \e^{-\ell}
\delta_{jb})  \kappa_{00} \sigma_{ij}^{\ell}$$ which acts on
$\Lambda$ by the scalar $$\sum_{1 \le i < j \le n } \sum_{\ell
=0}^{m-1}
 (\delta_{ia} - \e^{\ell} \delta_{ja} )(\delta_{ib} - \e^{-\ell}
\delta_{jb})  \kappa_{00} = \kappa_{00} \sum_{1 \le i < j \le n }
\sum_{\ell =0}^{m-1} (-\e^{\ell} \delta_{ja}\delta_{ib} -
\e^{-\ell} \delta_{ia} \delta_{jb})$$ which vanishes since
$\sum_{\ell=0}^{m-1} \e^{\ell} =0$.

Second, if $a=b$ then we must show that
$$1+
 \left[ \sum_{j=0}^{d-1}
(\kappa_{j+1} - \kappa_j) \sum_{r=0}^{d-1} \e^{prj} \right]  +
\sum_{1 \le i < j \le n } \sum_{\ell =0}^{m-1} (\delta_{ia} -
\e^{\ell} \delta_{ja} )(\delta_{ia} - \e^{-\ell} \delta_{ja})
\kappa_{00}$$ vanishes. So we must show that $1+ d \kappa_1 +
\kappa_{00} (m(n-a) + (a-1)m) =0$, which holds by the hypothesis.
\end{proof}

\subsection{A shift isomorphism}
We require the following theorem from \cite{BC}.

\begin{thrm}[Berest, Chalykh]\label{shift}
Let $(\kappa_{00}, \kappa_1, \ldots , \kappa_{d-1}) \in \C^d$ and
define $\kappa_{00}^{'} = \kappa_{00} +1 , \kappa_1^{'} = \kappa_1
+1$ and $\kappa_i^{'} = \kappa_i$ for all other $i$. Then there is
an isomorphism
$$\psi: e H_{\kappa^{'}} e \rightarrow e_\e H_\kappa e_\e$$
where $e$ is the symmetrising idempotent $e= \frac{1}{|W|} \sum_{w
\in W} w \in \C W$ and $e_\e = \frac{1}{|W|} \sum_{w \in W}
\det(w)  w$.
\end{thrm}
\begin{proof}
This is a minor modification of \cite[Theorem 5.8(2)]{BC}.
\end{proof}
\subsection{}
 Given $\kappa$ with $d \kappa_1 + m(n-1) \kappa_{00} =
-1-m(n-1)-d$, define $\kappa^{'}$ as in the statement of theorem
\ref{shift}. Then $d \kappa^{'}_1 + m(n-1) \kappa_{00}^{'} =   -1$
and hence by theorem \ref{onedim}, there is a one-dimensional
module $\Lambda$ for $H_{\kappa^{'}}$ on which $W$ acts trivially.
Hence, $e \Lambda = \Lambda$ becomes a $e H_{\kappa^{'}}
e$--module. Via the isomorphism of \ref{shift}, we may define a
$e_\e H_{\kappa} e_\e$--module $\Lambda^{\psi}$. Finally, we set
$$L = H_\kappa e_\e \otimes_{e_\e H_\kappa e_\e} \Lambda^{\psi},$$
an $H_\kappa$--module.

Then, because $H_{\kappa} e_\e$ is a finite $e_\e H_{\kappa}
e_\e$--module (this follows by considering the associated graded
modules, for instance), $L$ is finite-dimensional, and it then
follows from Section \ref{catoaltdef} that
$L$ belongs to the category $\cO$ of $H_{\kappa}$--modules.

Suppose we have chosen $\kappa$ generically. Then we may apply the
results of the previous section. In particular, Theorem
\ref{infdim} says that the only finite-dimensional simple object
in category $\cO$ is $L(\triv)$. Since $L$ has a composition
series with composition factors $L(\tau)$, $\tau \in \irrep(W)$,
we see that every composition factor of $L$ must be $L(\triv)$.
There is a functor $F: H_\kappa-\mathrm{mod} \rightarrow e_\e
H_\kappa e_\e-\mathrm{mod}$ defined by $FM = e_\e M$. This is an
exact functor and it takes a composition series of $L$ to a
composition series of $FL \cong \Lambda^{\psi}$. Hence, $e_\e
L(\triv) \neq 0$ and $L \cong L(\triv)$. This proves the following
lemma.
\begin{lem}\label{lequalsltriv}
If $\kappa_i$, $\kappa_{00}$ are chosen generically then the
$H_\kappa$--module $L$ is isomorphic to $L(\triv)$.
\end{lem}

\section{A quotient ring of the diagonal
coinvariants}\label{mainproof}
\subsection{} We follow the proof of \cite[Section 5]{Go} to obtain
the desired ring $S_W$ of Theorem \ref{mainthm}. Choose generic
$\kappa^{'}$ with $d \kappa_1^{'} + m(n-1)\kappa_{00}^{'}= -1$,
let $\kappa$ be defined as above, and define $L = H_{\kappa} e_\e
\otimes_{e_\e H_{\kappa} e_\e} \Lambda^\psi$ as above. Consider
the filtration on $H_{\kappa}$ with $\deg(\h) = \deg(\h^*) =1$ and
$\deg(W)=0$, and the associated graded module $\mathrm{gr} L$. As
in \cite{Go}, one obtains a surjection of $\C [\h \oplus
\h^*]*W$--modules:
$$\C[\h \oplus \h^*]*W e_\e \otimes_{\C [\h \oplus \h^*]^W}
\mathrm{gr} \Lambda^\psi \rightarrow \mathrm{gr} L.$$
By definition 
, $\C [\h \oplus \h^*]^W_+$ acts on $\mathrm{gr} \Lambda^\psi$ by
0, and hence $S_W:= \mathrm{gr} L \otimes \wedge^n \h$ is a
quotient ring of $D_W = \C[\h \oplus \h^*]/\langle \C[\h \oplus
\h^*]^W_+ \rangle$. We wish to determine the graded character of
$S_W$. To do so, we will determine the graded character of $L
\cong L(\triv)$. This requires a more delicate study of the
results of \cite{CE} that were used above.

\subsection{}
Recall from the proof of \ref{infdim} that we have the
$H_\kappa$--module $\tilde{Y}_c$ which is a finite-dimensional
quotient of the Dunkl representation $M(\triv)$ obtained by
factoring out an ideal $J$ of the polynomial ring $\C[\h] =
M(\triv)$, where $J$ is generated by a copy of $\h_q$ in degree
$m(n-1)+d+1$ consisting of singular vectors. Call this copy $U$.
Note that $U$ is an irreducible representation of $G(m,p,n)$ since
it follows from the definition of $\h_q$ (quoted in the proof of
\ref{infdim}) that $U \cong \h^*$ as $G(m,p,n)$--modules. Now we
are in a position to apply \cite[Theorem 2.3]{CE} with
$W=G(m,p,n)$. This says that we have a BGG-resolution
$$0 \leftarrow \tilde{Y}_c \leftarrow M(\triv) \leftarrow M(U)
\leftarrow \cdots \leftarrow M(\wedge^n U) \leftarrow 0.$$ Hence
in the Grothendieck group $K_0 (\cO)$, we have an equality
$$[\tilde{Y}_c] = \sum_{i=0}^n (-1)^i [M(\wedge^i \h^*)].$$
Now, by theorem \ref{infdim}, there is only one finite-dimensional
simple module in $\cO$, namely $L(\triv)$. Hence, $[\tilde{Y}_c] =
a \cdot [L(\triv)]$ for some $a \ge 1$. By \cite[Section 2.5
(32)]{DO}, there is an ordering on $\irrep(W)$ such that the
matrix with entries $[M(\tau): L(\sigma)]$ is unipotent upper
triangular. Therefore the classes $[M(\tau)]$ give a
$\mathbb{Z}$--basis for $K_0 (\cO)$ and there is a unique
expression $$[L(\triv)] = \sum_\tau c_\tau [M(\tau)]$$ with
$c_\tau \in \mathbb{Z}$. Hence, $a \cdot c_\triv = 1$ and so
$a=1$. It follows that $[\tilde{Y}_c] = [L(\triv)]$ and hence
$\tilde{Y}_c \cong L(\triv)$. Furthermore,
$$[L(\triv)] = \sum_{i=0}^n (-1)^i [M(\wedge^i \h^*)]$$ in $K_0 (\cO)$.

\subsection{}
 To prove Theorem
\ref{mainthm}, we will require the following lemma. Recall from
Section \ref{zdef} that the element $z \in \C W \subset
H_{\kappa}$ is defined by $z = \sum_{H \in
\mathcal{A}}\sum_{i=1}^{e_H-1} e_H k_{H,i} \e_{H,i}$. When
$W=G(m,p,n)$ we have
$$z = \sum_{i=1}^n \sum_{t=1}^{d-1} \kappa_t \sum_{j=0}^{d-1}
 \e^{ptj} s_i^{pj} + \kappa_{00} \sum_{i<j} \sum_{0 \le r \le m-1}
  (1 - \sigma_{ij}^{(r)} ).$$
\begin{lem}\label{zscalar}
$z$ acts on $\wedge^i \h^*$ by the scalar $i(d \kappa_1 + m(n-1)
\kappa_{00}).$
\end{lem}
\begin{proof}
Since $z$ is central, it acts on $\wedge^i \h^*$ by the scalar
$\chi_{\wedge^i \h^*}(z)/{n \choose i}$. Choosing a basis for $
\h^*$ for which $s_i^j$ is diagonal, we see that $\chi_{\wedge^i
\h^*}(s_i^j) = \e^j {n-1 \choose i-1} + {n-1 \choose i}$.
Similarly,  $\chi_{\wedge^i \h^*}(\sigma_{ij}^{(r)}) = {n-1
\choose i} - {n-1 \choose i-1}$. Substituting these values into
the expression for $z$ gives the result.
\end{proof}

\subsection{}
All of the properties of $S_W$ listed in Theorem \ref{mainthm}
apart from Theorem \ref{mainthm}(3) are immediate consequences of
the following theorem. Recall that $D_W$ is graded with
$\deg(\h)=-1$ and $\deg(\h^*)=1$ and that this grading is
$W$--stable. In general, if $M$ is a $\mathbb{Z}$--graded module
and $\chi_k$ is the character of the $k^{\mathrm{th}}$ graded
piece then the \emph{graded character} of $M$ is defined to be the
formal power series $\sum_i \chi_i t^i$.
\begin{thrm}\label{gradedchar}
The graded character of $\mathrm{gr} L = S_W \otimes \wedge^n
\h^*$ is
$$ w \mapsto t^{-n-m{n \choose 2}} \frac{ \det|_{\h^*}(1-t^{m(n-1)+d+1} w)}
{\det|_{\h^*}(1-tw)}$$
\end{thrm}
\begin{proof}
Recall the element $\mathsf{eu}_{\kappa} \in H_{\kappa}$ from
Section \ref{zdef}. By Lemma \ref{zscalar}, $z$ acts by 0 on the
trivial representation $\mathsf{triv}= \wedge^0 \h^*$ of $W$ and
hence $\mathsf{eu}_{\kappa}$ also acts by 0 on the trivial
representation $\mathsf{triv}$. Hence, the eigenvalue of
$\mathsf{eu}_{\kappa}$ on the subspace $\C[\h]_d \otimes
\mathsf{triv} \subset M(\mathsf{triv})$ is $d$.

By \cite[Theorem 4.2]{CE}, the representation
$\tilde{Y}_c=L(\mathsf{triv})$ of $H_{\kappa}$ is isomorphic as a
graded $G(m,1,n)$--module to $U_{m(n-1)+d+1}^{\otimes n}$ where
$U_{m(n-1)+d+1} = \C[u]/(u^{m(n-1)+d+1})$, regarded as a
representation of $\mathbb{Z}_m = \langle s_1 \rangle$ via $s_1
(u) = \e^{-1} u$, and where $S_n$ acts by permuting the factors of
the tensor product. Since $\{u^{i}| 0 \le i \le m(n-1)+d \}$ is a
basis of $U_{m(n-1)+d+1}$, we may define distinct basis elements
$a_i$ by $a_i := u^{m(i-1) +1}$, $1 \le i \le n$. Then the element
$v := \sum_{\sigma \in S_n} \mathrm{sgn}(\sigma) a_{\sigma(1)}
\otimes a_{\sigma(2)} \otimes \cdots \otimes a_{\sigma(n)}$
affords the representation $\wedge^n \h^*$ of $G(m,1,n)$, and lies
in degree $n+ \sum_{i=1}^n m(i-1) = n +m{n \choose 2}$. Hence,
$\mathsf{eu}_{\kappa} v =( n + m {n \choose 2}) v$. Note that $v$
also affords the representation $\wedge^n \h^*$ of $G(m,p,n)$ when
we consider $\tilde{Y}_c$ as a $W=G(m,p,n)$--module.

 We define ${\bf h} = \mathsf{eu}_{\kappa} - n -
m{n \choose 2} \in H_{\kappa}$.
 We now calculate the graded character of $L =
L(\mathsf{triv})$ with respect to the ${\bf h}$--eigenspaces. We
have shown above that in the Grothendieck group of $\cO$,
$$[L(\mathsf{triv})] = \sum_{i=0}^n (-1)^i [M(\wedge^i \h^*)].$$

Now, by Lemma \ref{zscalar}, $z$ acts by $-i(m(n-1)+d+1)$ on
$\wedge^i \h^*$, and hence, by Section \ref{zdef}, the lowest
eigenvalue of ${\bf h}$ on  $M(\wedge^i \h^*)$ is $i(m(n-1)+d+1) -
n -m{n \choose 2}$. Therefore, the graded character of $M(\wedge^i
\h^*)$ is $$t^{-n-m{n \choose 2}}\frac{ \chi_{\wedge^i \h^*}(w)
t^{i(m(n-1)+d+1)}}{\det|_{\h^*} (1-tw)}.$$ But $\det|_{\h^*} (1-
t^{m(n-1)+d+1}w) = \sum_i (-1)^i \chi_{\wedge^i \h^*} (w)
t^{i(m(n-1)+d+1)}$ (this follows readily from diagonalising $w$)
which gives the graded character of $L(\mathsf{triv})$ with
respect to the ${\bf h}$--eigenspaces as
$$w \mapsto t^{-n-m{n \choose 2}} \frac{ \det|_{\h^*}(1-t^{m(n-1)+d+1} w)}
{\det|_{\h^*}(1-tw)}$$ By the definition of the diagonal
coinvariant ring $D_W$, there is a unique copy of the trivial
representation in $D_W$, which lies in degree 0, and hence a
unique copy of $\wedge^n \h^*$ in $S_W \otimes \wedge^n \h^*$, and
hence a unique copy of $\wedge^n \h^*$ in $L(\mathsf{triv})$
(since $S_W \otimes \wedge^n \h^* = \mathrm{gr} L$, which is
isomorphic to $L$ as a $W$--module). The unique copy of $\wedge^n
\h^*$ in $L(\mathsf{triv})$ must be spanned by $v$. But ${\bf h} v
= 0$ and hence, by Lemma \ref{lequalsltriv}, ${\bf h}$ must act by
0 on the element $e_\e \otimes 1 \in L$ which affords the unique
copy of $\wedge^n \h^*$ in $L$. But because the grading induced by
${\bf h}$ on $L$ gives $e_\e \otimes 1$ degree 0 and $x \in \h^*$
degree 1, and $y \in \h$ degree $-1$, we see that $\mathrm{gr} L$
has the same graded character as $L(\mathsf{triv})$, which proves
the theorem.
\end{proof}

\subsection{Proof of Theorem \ref{mainthm} (3)} This is similar
 to a proof in \cite[Section 5]{Go}.
It is well-known (see for example \cite{kane}) that the ring of
coinvariants $A=\C[\h]/\C[\h]^W_+$ satisfies Poincar\'{e} duality.
Therefore the highest degree graded component of $A$, which lies
in degree $\sum_i (d_i-1)$ where the $d_i$ are the degrees of the
fundamental invariants of $W$, is an ideal of $A$ which is
contained in every nonzero ideal. This ideal is called the socle
of $A$.
In the case of $W=G(m,p,n)$ the degrees are $m, 2m, \ldots (n-1)m
, nd$ by \cite[Table 1, Table 2]{BMR}, so the socle lies in degree
$m {n \choose 2} +nd -n$. The image of $\mathbb{C}[\h]$ in $S_W$
corresponds to the subspace $\mathbb{C}[\h] e_{\e}\otimes
\Lambda^\psi$ of $L$. If $p\in \mathbb{C}[\h]^W_+e_{\e}$ then by
the definition of the shift isomorphism $\psi$ given in \cite{BC},
we see that $\psi(epe)= e_{\e}pe_{\e}$. It follows that in $L$ we
have:
$$p\otimes \Lambda^\psi = e_{\e}pe_{\e}\otimes \Lambda^\psi =
e_{\e}\otimes e_\e p e_\e \Lambda^\psi = e_\e \otimes epe \cdot
\Lambda = 0.$$
Thus the ideal generated by $\mathbb{C}[\h]_+^W$ annihilates
$e_{\e}\otimes \Lambda^\psi$. On the other hand, the quotient
$\mathbb{C}[\h]/\C[\h]^W_+$ contains a unique (up to scalar)
element of maximal degree $m {n \choose 2} +nd-n$, say $q$. The
space $\mathbb{C}q$ is the socle of $\mathbb{C}[\h]/\C[\h]^W_+$.
We claim $qe_{\e}\otimes \Lambda^\psi \neq 0$. By the PBW theorem,
any element of $H_\kappa$ can be written as a sum of terms of the
form $p_-wp_+$ where $p_-\in\mathbb{C}[\h^*]$,
$p_+\in\mathbb{C}[\h]$ and $w\in W$. Since $p_-$ and $w$ do not
increase degree, it would follow if $qe_{\e}\otimes \Lambda^\psi$
were zero, then $L$ could have no subspace in degree $m {n \choose
2}+nd-n$. But the Hilbert series of $L$ has highest order term
$t^{-n-m{n \choose 2} + mn(n-1) +nd} = t^{m{n \choose 2} +nd-n}$.
Thus $qe_{\e}\otimes \Lambda^\psi$ is non--zero and
$\mathbb{C}[\h]e_{\e}\otimes \Lambda^\psi$ is isomorphic to
$(\mathbb{C}[\h]/\C[\h]^W_+) e_\e \otimes \Lambda^\psi$. This
proves Theorem \ref{mainthm} (3). $\Box$

\section{Appendix: the case $W=G(m,p,2)$ with $p$ even}\label{fudge}
In this case, we have $d+1$ conjugacy classes of complex
reflections in $W$. In the notation of Section \ref{def} they are
$C_q:=\{s_i^{qp}|1 \le i \le n\}$, (for $1 \le q \le d-1$),
$C_{odd}:=\{ \sigma_{12}^{(\ell)} |  \text{$\ell$ odd}\}$ and
$C_{even}:=\{\sigma_{12}^{(\ell)}| \text{$\ell$ even} \}$. Thus,
the rational Cherednik algebra depends on $d+1$ complex
parameters: $\kappa_1, \kappa_2, \ldots, \kappa_{d-1}$
corresponding to $C_1$, $\kappa_{00}^{\rm{odd}}$ corresponding to
$C_{odd}$ and $\kappa_{00}^{\rm{even}}$ corresponding to
$C_{even}$. We choose these parameters so that
$\kappa_{00}^{\rm{odd}} = \kappa_{00}^{\rm{even}}=: \kappa_{00}$,
so our parameters become $(\kappa_{00}, \kappa_1, \ldots,
\kappa_{d-1})$. We must now verify that the rest of the proof of
Theorem \ref{mainthm} still works.

Using \cite[Table 1, Table 2]{BMR}, we see that in the definition
of the Hecke algebra $\mathcal{H}$ of Section \ref{hecke}, we get
the same braid relations and the same relation for $T_s$ as in
Section \ref{hecke} for the $p=2$ case, but the relations for the
$T_i$ are:
\begin{align*}
(T_{t_2^{'}}-1)(T_{t_2^{'}}+ e^{2 \pi i \kappa_{00}^{\rm even}})
&=0 \\
(T_{t_2} -1)(T_{t_2}+e^{2 \pi i \kappa_{00}^{\rm odd}}) &=0
\end{align*}
Thus, when $\kappa_{00}^{\rm odd}=\kappa_{00}^{\rm even}$, we get
the Hecke algebra of \cite[2.A]{GJ}. Furthermore, in \cite[Section
4.1]{CE}, it is assumed that
$\kappa_{00}^{\rm{odd}}=\kappa_{00}^{\rm{even}}$, so the
constructions of \cite{CE} are still valid. So the proofs of
Sections \ref{hecke} and \ref{hsimples} go through in the present
case. Hence, Theorem \ref{infdim} holds.

The only potential obstacle to completing the proof is the shift
isomorphism of Theorem \ref{shift}. But, under the shift
isomorphism $e H_{\kappa^{'}} e \rightarrow e_\e H_\kappa e_\e$
from \cite{BC}, the parameters $\kappa_{00}^{\rm{odd}}$ and
$\kappa_{00}^{\rm{even}}$ are both shifted by $1$, so we can
regard $\kappa_{00}$ as also being shifted by 1.
The rest
of the proof of Theorem \ref{mainthm} now goes through.

\bibliographystyle{alpha}
\bibliography{penguin}
\end{document}